\documentstyle[epsf, amssymb]{article}
\font\bbb=msbm10 scaled 1100

\newtheorem{thm}{Theorem}[section]
\newtheorem{lemma}[thm]{Lemma}
\newtheorem{cor}[thm]{Corollary}
\newtheorem{prop}[thm]{Proposition}

\newtheorem{dfn}[thm]{Definition}
\newtheorem{ex}[thm]{Example}

\newtheorem{rem}[thm]{Remark}

\newcommand{\eps}{\epsilon}

\newcommand{\la}{\lambda}			
\newcommand{\La}{\Lambda}
\newcommand{\F}{{\cal F}}
\newcommand{\del}{\partial}

\newcommand{\ra}{\rightarrow}

\newcommand{\hra}{\hookrightarrow}
\newcommand{\real}{\mbox{\bbb R}} 	
\newcommand{\zed}{\mbox{\bbb Z}}	
\newcommand{\inv}{^{-1}}

\newcommand{\grad}{\nabla}

\newcommand{\norm}[1]{\|#1\|}

\newcommand{\rest}[2]{\left. #1\right\vert_{#2}}		
\newcommand{\be}{\begin{equation}}
\newcommand{\ee}{\end{equation}}
\newcommand{\bea}{\begin{eqnarray}}
\newcommand{\eea}{\end{eqnarray}}
\newcommand{\bmini}{\footnotesize\begin{center}\begin{minipage}{5.5in}}
\newcommand{\emini}{\end{minipage}\end{center}\normalsize}
\newcommand{\pf}{{\em Proof: }}

\newcommand{\qed}{\hfill$\Box$}

\renewcommand{\L}{{\cal L}}

\newcommand{\W}{{\cal W}}

\newcommand{\bfv}{\mbox{{\bf v}}}

\newcommand{\eg}{{\em e.g.}}
\newcommand{\ie}{{\em i.e.}}
\newcommand{\cf}{{\em cf. }}
\begin{document}

\begin{center}
\Large 
{\bf   Gradient flows within plane fields}
\normalsize
\vspace{0.1in}

\begin{tabular}{cc}
JOHN ETNYRE & ROBERT GHRIST \\ 
Department of Mathematics & School of Mathematics \\
Stanford University & Georgia Institute of Technology \\
Stanford, CA 94305 & Atlanta, GA 30332-0160 
\end{tabular}
%
\end{center}


\begin{abstract}
We consider the dynamics of vector fields on three-manifolds which are
constrained to lie within a plane field, such as occurs in 
nonholonomic dynamics. On compact manifolds, such vector 
fields force dynamics beyond that of a gradient flow, except in cases 
where the underlying manifold is topologically simple 
(\ie, a graph-manifold). 
Furthermore, there are strong restrictions on the types of 
gradient flows realized within plane fields: such flows lie 
on the boundary of the space of nonsingular Morse-Smale 
flows. This relationship translates 
to knot-theoretic obstructions for the link of singularities
in the flow. In the case of an integrable plane field, the 
restrictions are even finer, forcing taut foliations on 
surface bundles.
The situation is completely different in the case of contact plane 
fields, however: it is easy to realize gradient fields within 
overtwisted contact structures (the nonintegrable 
analogue of a foliation with Reeb components).

\vspace{0.1in}

\vspace{0.1in}
\noindent{\bf keywords:} Gradient flows, foliations, contact structures, 
		round handles, knots.

\end{abstract}

\section{Introduction}
\label{sec_1}

Let $M$ denote a compact oriented three-manifold. A {\em plane field} 
$\eta$ on $M$ is a subbundle of the tangent bundle
$TM$ which associates smoothly to each point $p\in M$ a 
two-dimensional subspace $\eta(p)\subset T_pM$. Unlike line fields, 
a plane field cannot always be integrated to yield a two-dimensional 
foliation $\F$.
A plane field is said to be {\em integrable} if it can be ``patched
together'' to yield a foliation whose leaves are tangent to the plane field at
each point. Certainly, such plane fields have strong topological and
geometric properties. On the other hand, the 
case where the plane field $\eta$ is nowhere integrable can be equally 
important. A maximally nonintegrable (in the sense of Frobenius --- 
see \S\ref{sec_5}) plane field on an odd-dimensional
manifold is a {\em contact structure}. Seen as an ``anti-foliation'',
contact structures are rich in geometric and topological properties
which of late have become quite important in understanding
the topology of three-manifolds and the symplectic
geometry of four-manifolds.

Let $X$ be a vector field on $M$. The dynamics of $X$ are often  
related to global properties of $M$. If we further specify that $X$ 
is tangent to a plane field $\eta$ --- that is, 
$X(p)\in\eta(p)$ for all $p\in M$ 
--- then we might expect stronger relationships. 
We will consider the ways in which the topology and geometry of 
a plane field $\eta$ are coupled to the dynamics of vector fields 
contained in $\eta$. The general principal at work here as elsewhere
is that {\em simple dynamics implicate simple topological objects in 
dimension three.} We will reassert this by examining the gradient 
flows within plane fields.

The examination and classification of gradient flows  
has been ubiquitous in the study of manifolds: \eg, the h-cobordism 
theorem and the resolution of the high-dimensional 
Poincar\'e Conjecture. 
This paper will add to the typical scenario the constraint of lying 
within a plane field. Atypical restrictions on the 
dynamics and on the underlying manifold are born out of this.

We note that the problem of understanding gradient fields 
constrained to lie within plane fields is by no means unnatural. 
The study of mechanical systems with nonholonomic constraints
is precisely the study of flows constrained to lie within a 
nowhere integrable distribution (\ie, in odd dimensions, a 
contact structure). For example, gradient flows for mechanical systems have 
been used successfully in the control of robotic systems (see, \eg, 
\cite{KR90}): to maneuver a robot from points $A$ to $B$ through 
a physical space replete with obstacles, one establishes a gradient flow on 
a suitable configuration space with $B$ as a sink, with $A$ in
the basin of attraction for $B$, and with infinite
walls along the obstacles. In this paper, we show that the 
nonholonomic version of this procedure possesses potentially 
difficult topological obstructions. 

The paper is organized as follows: the remainder of this section 
provides a brief sketch of the requisite theory from the dynamical 
systems approach to flows. In 
\S\ref{sec_2}, we commence our investigation of plane field 
flows by examining local and global properties of fixed points: 
fixed points will not be isolated, but must (on an open dense 
subset of $C^r$ vector fields tangent to $\eta$, $r\geq 1$)
rather appear in {\em links}, or embedded closed curves. 
This culminates in a classification of gradient flows on three-manifolds 
which can lie within a plane field in \S\ref{sec_3}. 
The existence of such flows is 
equivalent to the existence of a certain type of 
round handle decomposition for the manifold 
(see Definition \ref{def_RH}). Surprisingly, this same 
restriction appears when considering energy surfaces for 
(Bott-) integrable Hamiltonian flows \cite{CMAN}.

\noindent{\bf Theorem:} {\em
Let $M$ be a compact 3-manifold outfitted with a plane field
$\eta$. If $X$ is a nondegenerate\footnote{See 
Definition~\ref{def_Nondegenerate}.} gradient field tangent to $\eta$, 
then $X$ lies in the boundary of the space of nonsingular 
Morse-Smale flows on $M$. Furthermore, the set of fixed 
points for $X$ forms the cores of an essential round handle 
decomposition for $M$.
}

This leads to the corollary (a stronger form of which is proved 
in \S\ref{sec_3}):

\noindent{\bf Corollary:} {\em 
Non-gradient dynamics is a generic (residual) property in the class of 
$C^r$ ($r\geq4$) vector fields tangent to a fixed $C^r$ 
plane field on a closed hyperbolic three-manifold.
}

In \S\ref{sec_4}, we consider the manifestation of these restrictions
on a knot-theoretic level for the particular case of the 3-sphere. 

\noindent{\bf Theorem:} {\em
For $X$ a nondegenerate gradient plane field flow on $S^3$, 
each connected component of the fixed point set of $X$ is a knot 
whose knot type is among the class generated from the unknot 
by the operations of iterated cabling and connected sum.
}

We proceed with remarks on two cases in which the plane field carries
additional geometric structure: first, the case of an everywhere integrable 
plane field, \ie, a foliation; and second, the case of a maximally 
nonintegrable plane field, \ie, a contact structure. 
The property of carrying a gradient flow in a foliation forces the 
foliation to be taut; hence, there are no 
(nondegenerate) gradient flows within a foliation on $S^3$. 
More generally, we have the following restrictions on the underlying
three-manifold:

\noindent{\bf Theorem:} {\em
A closed orientable three-manifold containing a nondegenerate gradient 
field within a $C^r$ ($r\geq2$) codimension-one foliation must be a surface 
bundle over $S^1$ with periodic (or reducibly periodic) monodromy map.
}

The corresponding restrictions do not hold for the contact case. We 
demonstrate that gradient fields can always reside within the analogue 
of a non-taut foliation: an {\em overtwisted} contact structure.
We close with two questions on the higher dimensional 
versions of the results of this paper.

\subsection{The dynamics of flows}
\label{sec_1.1}

Ostensibly, flows within a plane field would appear to be a relatively
restricted class of objects. However, the dynamics of such flows can
exhibit behaviors which range from strictly two-dimensional
dynamics (as when the plane field yields a foliation by compact
leaves) to fully three-dimensional phenomena (\eg, an Anosov
flow, which is tangent to a pair of transverse integrable plane 
fields). In \S\ref{sec_2}, we show that near a fixed point of 
a plane field flow, the dynamics are locally ``stacked'' planar
dynamics. In contrast, it is a simple exercise in homotopy 
theory that every nonsingular flow on $S^3$ (or any integral  
homology 3-sphere) lies within a plane field. 

A few definitions are important for the dynamical systems 
theory used in this paper.
The most important aspect of a flow with respect to its 
geometry and dynamics is the notion of hyperbolicity.
Recall that an invariant set $\La\subset M$ of a flow
$\phi^t$ is {\em hyperbolic}
if the tangent bundle $\rest{TM}{\La}$ 
has a continuous $\phi^t$-invariant splitting 
into $E^\phi\oplus E^s\oplus E^u$, 
where $E^\phi$ is tangent to the flow direction, and $D\phi^t$  
uniformly contracts and expands along $E^s$ and $E^u$ 
respectively: \ie,
\be
\begin{array}{ll}
\norm{D\phi^t(\bfv^s)} \leq Ce^{-\la t}\norm{\bfv^s} & 
	\mbox{for } \bfv^s\in E^s \\
\norm{D\phi^{-t}(\bfv^u)} \leq Ce^{-\la t}\norm{\bfv^u} & 
	\mbox{for } \bfv^u\in E^u \\
\end{array}, t>0,
\ee
for some $C\geq 1$ and $\la > 0$.
A flow $\phi^t$ which is hyperbolic on all of $M$ is 
called an {\em Anosov flow.}

The existence of hyperbolic invariant sets greatly simplifies 
the analysis of the dynamics. The principal tool available  
is the Stable Manifold Theorem \cite{HPS70}, which states 
that for a hyperbolic invariant set, the distributions 
$E^s$ and $E^u$ are in fact tangent to 
global {\em stable} and {\em unstable manifolds}: 
manifolds, all of whose points have the 
same backwards and forwards (resp.) asymptotic behavior. 
See any of the standard texts (\eg, \cite{GH83}) 
for further information and examples.

\section{Fixed points}
\label{sec_2}

In analyzing the dynamics and topology of a flow, one examines
dynamical $n$-skeleta of increasing dimension: 
first the fixed points, then periodic and connecting orbits,
lastly higher-dimensional invariant manifolds and attractors. 
This section concerns the typical distribution of fixed points 
for plane field flows.

\begin{lemma}\label{lem_Coords}
Given $\eta$ a $C^r$ plane field on $M^3$ and $p\in M$ there exists a
neighborhood $U\cong\real^3$ of 
$p$ along with local coordinates $(x,y,z)$ on $U$ such that 
$\eta=\mbox{ker}(\alpha)$, where $\alpha$ is a one-form given by
\be \label{eq_normal}
	\alpha = dz + g(x,y,z)dy, 
\ee
for some function $g$ which vanishes at the origin. 
The space $\Gamma^r(\rest{\eta}{U})$ of $C^r$ sections of $\eta$  
on $U$ is isomorphic to $C^r(\real,C^r(\real^2,\real^2))$, 
the space of $C^r$ arcs of $C^r$ planar vector fields. 
\end{lemma}
\pf
That $\alpha$ exists is easy to derive (and is stated in 
\cite{ET97}): choose coordinates $(x,y,z)$ so that $\del/\del z$ is 
transverse to $\eta$ on $U$. Then, after rescaling, $\eta$ is the kernel of 
$dz+f(x,y,z)dx + g(x,y,z)dy$. By a change of variables, one can 
eliminate $f$ and remove constant terms in $g$.

Parameterize $U$ as $\{\real^2\times\{z\} : z\in\real\}$. 
Given any 1-parameter family of functions $F_z:\real^2\ra\real^2$, 
there is a well-defined vector field on $U$ given by
\be\label{eq_Coords}
\begin{array}{l}
	\dot{x} = f_1(x,y,z) \\
	\dot{y} = f_2(x,y,z) \\
	\dot{z} = -g(x,y,z)f_2(x,y,z) 
\end{array}
\; \mbox{ where } F_z(x,y)=\left( f_1(x,y,z) , f_2(x,y,z) \right) ,
\ee 
which lies within $\eta$ by Equation~\ref{eq_normal}. Similarly, any 
vector field on $U$ contained in $\eta$ induces a 1-parameter family 
of planar vector fields $F_z:\real^2\ra\real^2$ by inverting the above 
procedure. Since $\del/\del z$ is always transverse to $\eta$, 
zeros of $F_z$ correspond precisely with zeros of 
the induced vector field in $\eta$. Note finally that the correspondence 
is natural with respect to the $C^r$-topology 
(nearby families of planar vector fields induce nearby plane field 
flows and vice versa).
\qed

\begin{prop}\label{prop_Link}
Let $\eta$ be a $C^r$ ($r\geq 1$) plane distribution on 
$M$ a compact 3-manifold, 
and let $\Gamma(\eta)$ denote the space of $C^r$ sections
of $\eta$. Then on an open dense subset of $\Gamma(\eta)$, 
the fixed point set is a smooth finite link of embedded circles.
\end{prop}
\pf
From the Transversality Theorem (see \cite[p. 74]{Hir76})
we know there is an open dense subset of sections of $\eta$
which are transverse to the zero section.  The proposition 
clearly follows.
\qed

\begin{cor}\label{cor_NoHypFP}
Let $X$ be any vector field on $M^3$ contained in the distribution 
$\eta$. Then any fixed point of $X$ is nonhyperbolic.
\end{cor}
\pf
Hyperbolic fixed points are isolated and persist in 
$C^1$-neighborhoods of vector fields; hence, they cannot be 
perturbed to yield circles of fixed points.
\qed

To analyze the dynamics near a curve of singularities, we show that
for all but finitely many points, the dynamics are 
{\em transversally hyperbolic}; \ie, after ignoring the 
nonhyperbolic direction along the curve, the flow 
is hyperbolic along the tangent plane transverse to the curve. 
We then turn to classify the (codimension-1) bifurcations 
in the transverse behavior along a curve of singularities.

\begin{prop}\label{prop_TransHyp}
Let $X\subset\eta$ be a $C^r$ ($r\geq2$) section of a $C^r$ plane 
field $\eta$. Then on a residual set of such vector fields,  
${\mbox{Fix}}(X)$ is a link $L$ which is transversally hyperbolic with 
respect to all but finitely many $p\in L$.
\end{prop}
\pf
By a standard argument (see \cite[p. 74]{Hir76}) is suffices to show
that there is an open cover $\{U_i\}$ of $M$ for which there is a 
residual set of sections of $\eta\vert_{U_i}$ with the desired property.
Cover each $p\in M$ by a chart as in Lemma~\ref{lem_Coords}. On 
each chart, consider the map from $\real^3\ra\real^2$ induced by 
a section of $\eta$. Extend this to a map into the 1-jet space 
$J^1(\real^3,\real^2)$ to capture information about the linearization 
of the flow. One may easily find a codimension three stratified subset $S$ of
$J^1(\real^3,\real^2)$ on which a section will both vanish and be 
transversally nonhyperbolic.  Thus by the Jet Transversality Theorem for
$C^2$ maps we obtain a residual subset of sections of $\eta$ whose 
1-jets transversally intersect $S$ at isolated points (which clearly 
must lie on $L$).
\qed

\begin{cor}\label{cor_QuadTan}
Under the hypotheses of Proposition~\ref{prop_TransHyp}, 
the singular link $L$ is transverse to $\eta$ at all but a finite number
of points.
\end{cor}
\pf
If the curve of singularities $S$ is tangent to the plane field 
$\eta$ at a point $p$, then $p$ is not transversally hyperbolic since the 
eigenvalue whose eigenvector points in the direction transverse 
to $\eta$ is zero (the vector field can have 
no component in the direction transverse to $\eta$). 
\qed

It is now a simple matter to classify the points at which the vector
field is not transversally hyperbolic to the equilibria. Thanks 
to Lemma~\ref{lem_Coords}, this analysis reduces simply 
to bifurcation theory of fixed points in planar vector fields.
In particular, there are precisely two ways in which a (generic) 
$X\subset\eta$ can fail to be transversally hyperbolic at a point.

Given any singular point $p\in S$, the transverse dynamics is 
characterized by the pair of transverse eigenvalues for the 
linearized flow: $\la^x$ and $\la^y$. Transverse hyperbolicity 
fails if and only if one or both of these eigenvalues has zero real part.
Generically, this can occur in two distinct ways. First, $\la^x$ and $\la^y$ 
may be both real, and one of them goes transversally through zero: 
this is a {\em saddle-node bifurcation}. Second, 
$\la^x$ and $\la^y$ may be a complex conjugate pair of eigenvalues
which together pass through the imaginary axis transversally: this is 
a {\em Hopf bifurcation}. Again, these names correspond
with analogous bifurcations of fixed points in planar vector fields.

\begin{figure}[htb]
\begin{center}
	\epsfxsize=3.0in\leavevmode\epsfbox{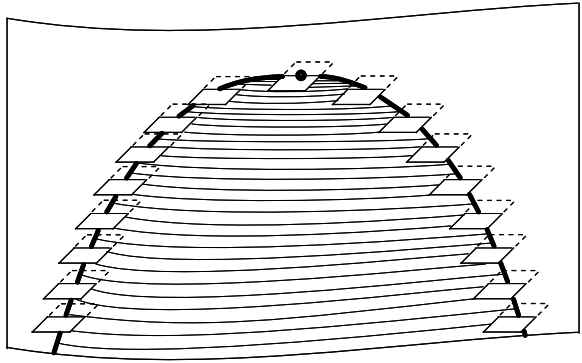}
\end{center}
\caption{A saddle-node bifurcation of singularities in a plane field
	flow.} 
\label{fig_SN}
\end{figure}

\begin{prop}\label{prop_SN}
In the unfolding of a $C^r$-generic ($r\geq2$) saddle node bifurcation 
on a curve of fixed points in a plane field flow, there is a 
quadratic tangency between the plane field and the fixed 
point curve, along with a one-parameter family of 
heteroclinic connections between fixed points limiting onto the
bifurcation point, as in Figure~\ref{fig_SN}.
\end{prop} 
\pf
As per Lemma~\ref{lem_Coords}, choose a coordinate system 
$(x,y,z)$ on a neighborhood of the bifurcation point $p$ 
so that $\del/\del z$ is everywhere transverse to the plane field 
$\eta$. It is also clearly possible (via the Stable Manifold 
Theorem) to choose coordinates so that 
the $x$-direction corresponds to the eigenvector for the 
transversally hyperbolic eigenvalue $\la^x$.

By Lemma~\ref{lem_Coords}, the unfolding of this 
codimension-1 fixed point in a plane field flow corresponds to 
the codimension-1 unfolding of a generic fixed point in a planar 
vector field having one hyperbolic eigenvalue and one 
eigenvalue with zero real part. The unfolding of the planar
saddle-node is conjugate to the system \cite{GH83}
\be
F_z : \begin{array}{l}
	\dot{x} = \la^x x \\
	\dot{y} = z - a y^2
\end{array} ,
\ee
for some $a\neq 0$, which, under Equation~\ref{eq_Coords}, 
corresponds to the vector field within $\eta$
\be
\begin{array}{l}
\dot{x} = \la^x x \\
\dot{y} = z  - a y^2 \\
\dot{z} = -g(x,y,z) (z - a y^2)
\end{array} .
\ee
The curve of fixed points is thus a parabola tangent to $\eta$
at the bifurcation point. 

To show the existence of a family of heteroclinic curves from 
one branch of the parabola to the next, note that the 
planar vector fields $F_z$ have precisely this 1-parameter 
family of orbits. Upon ``suspending'' to obtain a vector field
within $\eta$, the orbits remain, since the expression for 
$\dot{z}=-g(x,y,z)(z-a y^2)$ vanishes at $(0,0,0)$; hence,  
$\dot{z}$ is bounded near zero in a neighborhood of the bifurcation 
value and the integral curves within the invariant 
plane $\dot{x}=0$ must connect.
\qed

\begin{figure}[htb]
\begin{center}
	\epsfxsize=3.5in\leavevmode\epsfbox{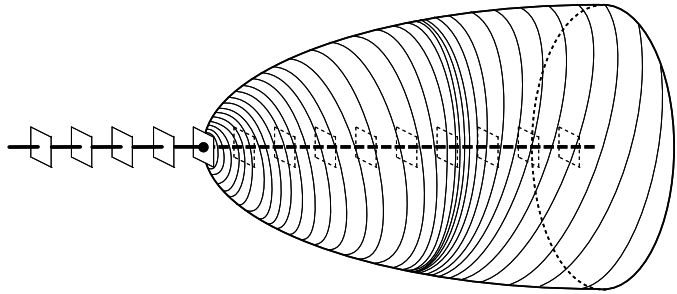}
\end{center}
\caption{A Hopf bifurcation of singularities in a plane field 
	flow.} 
\label{fig_Hopf}
\end{figure}

\begin{prop}\label{prop_Hopf}
In the unfolding of a $C^r$-generic ($r\geq4$) 
codimension-one Hopf bifurcation on a 
curve of fixed points in a plane field flow, there is an 
invariant attracting or repelling paraboloid which opens 
along the curve of fixed points as in Figure~\ref{fig_Hopf}. 
\end{prop}
\pf
Since only the real portion of the transverse eigenvalues
vanish, the curve of fixed points is transverse to the plane field
in a neighborhood of the bifurcation point $p$.
Hence, choose coordinates as per Lemma~\ref{lem_Coords} 
such that the curve of fixed points is the $z$-axis and the 
bifurcation point is at $(0,0,0)$. Again, by Lemma~\ref{lem_Coords}, 
this bifurcation in a plane field flow 
corresponds precisely to the codimension-one Hopf bifurcation 
of planar vector fields, conjugate to the truncated normal form 
\cite{GH83}
\be
F_z : \begin{array}{l}
\dot{r} = z r + a r^3 \\
\dot{\theta} = \omega
\end{array} ,
\ee
where we have transformed $(x,y)$ to polar coordinates   
and the constant $a$ is (in the codimension-1 scenario) nonzero. 
Solving this equation for $\dot{r}=0$ yields the paraboloid 
$r = \sqrt{-{z}/{a}}$, 
which is either attracting or repelling, depending on the sign 
of the coefficient $a$. By translating Lemma~\ref{lem_Coords} 
into polar coordinates, it follows that $\dot{z}$ is of
order $r^2$, which is less than $\dot{r}$; hence, 
adding the dynamics in the $z$-component
affects neither the existence nor the attracting/repelling nature
of the invariant paraboloid; however, unlike the planar case, 
the paraboloid is not necessarily fibered with closed curves. In 
general, orbits will spiral about the paraboloid.
\qed

\begin{rem}{\em \label{rem_Bifn}
We note that saddle-node or Hopf bifurcations must occur in 
pairs, since the fixed point curves are circles and the index at 
a bifurcation changes. However, in 
the case where there are no saddle-node 
or Hopf bifurcations along the singular curve, the flow is 
everywhere transversally hyperbolic, and the index of 
the fixed points (source, saddle, or sink) is constant along
the curve. 
}\end{rem}

We conclude with the definition of a nondegenerate vector field 
tangent to a plane field, and prove that such vector fields are
generic.

\begin{dfn}\label{def_Nondegenerate}{\em
A {\em nondegenerate} section of a plane field $\eta$
is a vector field $X\subset\eta$ whose fixed point set is a link having 
transversally hyperbolic dynamics at all but a finite number of 
points, at which the degeneracies are codimension one. 
}\end{dfn}

\begin{prop}
Nondegenerate fields are generic (residual in the $C^r$ topology
$r\geq4$) within the space of sections to a $C^r$ plane field $\eta$. 
\end{prop}
\pf
We simply repeat the argument in the proof of Proposition~\ref{prop_TransHyp} using
Propositions \ref{prop_SN} and \ref{prop_Hopf}.
\qed

\section{Round handles and gradients}
\label{sec_3}

Let $X$ be a nondegenerate vector field contained in the plane field $\eta$. 
The goal of the remaining sections is to understand restrictions 
on the topology of 3-manifolds supporting plane field flows which 
are forced by prescribed dynamics. A well-known example of this
occurs in the case of Anosov flows: certain three-manifolds are prohibited 
from carrying Anosov dynamics. In contrast, we examine obstructions 
associated to the simplest kinds of dynamics: gradient plane field flows. 
We show that only certain topologically ``simple'' manifolds support 
such dynamics. This will lead us to further knot-theoretic obstructions 
based on the singular links in a plane field flow. An old theme is 
played out: when the dynamics of $X$ are simple, the links associated 
to it are simple.

\begin{lemma}\label{lem_Nonsense}
Let $M$ denote an oriented Riemannian 3-manifold and $X=-\grad\Psi$ a 
$C^r$ ($r\geq2$) gradient vector field which lies within a $C^r$ 
plane field $\eta$ on $M$.
Then $\Psi$ is constant on each connected component of $\mbox{Fix}(X)$, the
fixed point set of $X$. Furthermore, if $c$ is a regular 
value of $\Psi$, then $\Psi\inv(c)$ is a disjoint union of tori
transverse to both $X$ and $\eta$.
\end{lemma}
\pf
Each component of $\mbox{Fix}(X)$ is a compact connected set 
of critical points for $\Psi$, whose image under $\Psi$ is a compact 
connected subset of $\real$ having measure zero, by the Morse-Sard
Theorem. For $c$ regular, $\Psi\inv(c)$ is a disjoint union of smooth 
surfaces, and $X$ is transverse to each component since $X$ is a
gradient field. Hence, the plane field $\eta$ is everywhere transverse
to $\Psi\inv(c)$ and the resulting line field given by the 
intersection of $\eta$ and the tangent planes to $\Psi\inv(c)$ in 
$\rest{TM}{\Psi\inv(c)}$ is nonsingular. Thus, the Euler 
characteristic of each component of $\Psi\inv(c)$ is zero. 
The transverse vector field $X$ gives an orientation to 
the surface, which excludes from consideration the Klein bottle.
\qed

Grayson and Pugh \cite{GP93} give examples of $C^\infty$ 
functions on $\real^3$ 
whose critical points consist of a smooth link, yet for which the level
sets are usually not tori: see Remark \ref{rem_GP}.

The above mentioned restrictions on gradient plane fields 
translate into very precise conditions on the topology of the 
underlying three-manifold. The fact that the manifold consists of 
a finite number of thick tori $T^2\times[0,1]$ glued together 
in ways prescribed by $\Psi$ implies that the manifold 
can be decomposed into solid tori in a canonical fashion: this 
phenomenon was identified and analyzed by Asimov and Morgan
in the 1970's \cite{Asi75,Mor78} in a completely different context.

\begin{dfn}{\em \label{def_RH} 
A {\em round handle} (or RH) in dimension three is a solid torus
$H=D^2\times S^1$ with a specified index and exit set
$E\subset T^2=\del(D^2\times S^1)$ as follows:
\begin{description}
\item[index 0:] $E = \emptyset$.
\item[index 1:] $E$ is either (1) a pair of disjoint annuli on the boundary
	torus, each of which wraps once longitudinally; or 
	(2) a single annulus which wraps twice longitudinally.
\item[index 2:] $E=T^2$. 
\end{description}
}\end{dfn}

\begin{dfn}{\em \label{def_RHD}
A {\em round handle decomposition} (or RHD) for a manifold $M$ is a 
finite sequence of submanifolds
\be \emptyset = M_0\subset M_1\subset\cdots M_n=M ,\ee
where $M_{i+1}$ is formed by adjoining a round handle 
to $\del M_i$ along the exit set $E_{i+1}$ of the round handle.
The handles are added in order of increasing index. 
}\end{dfn}

Asimov and Morgan \cite{Asi75,Mor78} 
used round handles to classify nonsingular
{\em Morse-Smale} vector fields: that is, vector fields whose
recurrent sets consist entirely of a finite number of hyperbolic
closed orbits with transversally intersecting invariant
manifolds. 

\begin{thm}\label{thm_RH}
Let $\eta$ denote a $C^r$ ($r\geq2$) plane field on $M^3$ (compact) with $X\subset\eta$ 
a $C^r$ nondegenerate gradient vector field. 
Then the set of fixed points for $X$ forms the cores of
a round handle decomposition for $M$. Furthermore, the 
indices of the fixed points correspond to the indices of 
the round handles, and $X$ is transverse to $\del M_i$ for all $i$. 
\end{thm}
\pf
Let $L$ denote the set of fixed points for $X=-\grad\Psi$: 
this is an embedded link. 
We first show that every fixed point is transversally hyperbolic. 
From Remark~\ref{rem_Bifn}, the only non-transversally 
hyperbolic points must occur as Hopf bifurcations or 
saddle-node bifurcations. Hopf bifurcations are associated to 
complex transverse eigenvalues, which cannot exist in a 
gradient flow. Similarly, a 
saddle-node bifurcation introduces a one-parameter family
of heteroclinic connections as in Figure~\ref{fig_SN}. This 
also cannot occur in a gradient flow, since by 
Lemma~\ref{lem_Nonsense} we have the function $\Psi$ constant 
on the curve of fixed points. The orbits of the flow 
which necessarily connect one side to the other cannot
be obtained by flowing down a gradient. Hence, 
each singular curve is transversally hyperbolic with constant index.

Choose $N$ a small tubular neighborhood of $L$ in $M$ and let $f$ 
denote a bump function in $N$ which evaluates to 1 on $L$
and is zero outside of $N$. Orient the link $L$ and perturb
$X$ to the new vector field $X+\eps f\frac{\del}{\del z}$, where 
$\frac{\del}{\del z}$ denotes the unit tangent vector along $L$. 
This yields a nonsingular flow which has $L$ as a set of 
hyperbolic closed orbits and no other recurrence. After a 
slight perturbation to remove any nontransverse intersections
of stable and unstable manifolds to $L$, this vector field 
is a nonsingular Morse-Smale field with periodic orbit link $L$.
The work of Morgan \cite{Mor78} then implies that $L$ forms
the cores of a round handle decomposition for $M$, where the 
index of each handle corresponds to the transverse index of 
the curve of fixed points (source, saddle, or sink). 
In \cite{Mor78} it is moreover shown that the nonsingular 
Morse-Smale vector field is transverse to each $\del M_i$; since 
the neighborhood $N$ is very small, this transversality 
remains in effect for $X$.
\qed

\begin{cor}
Gradient flows on plane fields in three-manifolds 
lie on the boundary of the space of nonsingular Morse-Smale fields.
\end{cor}
\pf
In the proof of Theorem~\ref{thm_RH}, let $\eps\ra 0$. This gives a 
one-parameter family of nonsingular Morse-Smale flows which 
converges to the gradient plane-field flow.  
\qed

\begin{cor}\label{cor_Non-Graph}
Non-gradient dynamics is a generic condition in the 
space of plane field flows on an irreducible non-graph three-manifold 
(\eg, a hyperbolic 3-manifold).
\end{cor}
\pf
By the work of Morgan \cite{Mor78}, round handle decompositions
of irreducible three-manifolds exist only for the class of graph-manifolds. 
\qed

Recall that a graph manifold is a three-manifold given by gluing together
Seifert-fibered spaces along essential torus boundaries. 
Examples include $S^3$, lens spaces, and manifolds with many 
$S^2\times S^1$ connected summands. 
The property of being composed of Seifert-fibered pieces (\ie, 
a graph manifold) is relatively rare among three-manifolds, 
the ``typical'' irreducible three-manifold being composed of 
hyperbolic pieces.

\begin{rem}{\em
We may push Theorem~\ref{thm_RH} a bit further. Let $\phi^t$ 
be a plane field flow whose chain-recurrent set consists entirely 
of transversally hyperbolic curves of fixed points and a finite
set of hyperbolic periodic orbits (note that hyperbolic periodic 
orbits can easily live within plane fields, even within nowhere 
integrable plane fields). This situation is, 
after the class of gradient flows, 
the next simplest scenario dynamically. Then, by the same proof, the 
connected components of the entire chain-recurrent set must form the 
cores of a round-handle decomposition. Hence, the additional 
dynamics forced upon plane field flows in a non-graph 
manifold is something other than hyperbolic periodic orbits.
}\end{rem}

\section{The link of singularities}
\label{sec_4}

\begin{figure}[htb]
\begin{center}
	\epsfxsize=4.5in\leavevmode\epsfbox{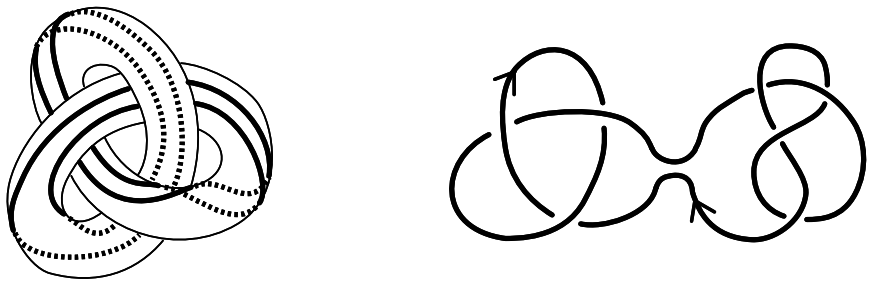}
\end{center}
\caption{Operations to generate zero-entropy knots: (left) cabling;
(right) connected sum.} 
\label{fig_Knots}
\end{figure}

We have shown that fixed points of plane field flows appear in 
links. The natural question is which links can 
arise as the singular points, and what dependence is there upon
the dynamics of the plane field flow. For nondegenerate gradient
fields, it is an immediate corollary of Theorem~\ref{thm_RH} that
the singular link is a collection of fibers in the 
Seifert-fibered portions of a graph manifold. We can be 
more specific, however, in the special case of $S^3$. We recall two
standard operations for transforming simple knots into more 
complex knots: see Figure~\ref{fig_Knots} for an illustration.

\begin{dfn}{\em
Let $K$ be a knot in $S^3$. Then the knot 
$K'$ is said to be a $(p,q)$-{\em cable} of
$K$ if $K'$ lives on the boundary of a tubular neighborhood of $K$,
wrapping about the longitude (along $K$) $p$-times and about the 
meridian (around $K$) $q$-times. 
Let $K$ and $J$ be a pair of knots in $S^3$. Then the {\em connected sum},
denoted $K\# J$, is defined to be the knot obtained by removing 
from each a small arc and identifying the endpoints along a 
band as in Figure~\ref{fig_Knots}.
}\end{dfn}

\begin{dfn}{\em
The {\em zero-entropy knots} are the collection of knots 
generated from the unknot by the operations of cabling and 
connected sum; \ie, it is the minimal class of knots closed under
these operations and containing the unknot. 
}\end{dfn}

Zero-entropy knots are relatively rare among all knots: 
\eg, none of the {\em hyperbolic knots} (such as the 
figure-eight knot, whose complement has a hyperbolic structure) 
are zero-entropy. The title stems from the often-discovered 
fact (see \cite{G97CSF} for history) 
that such knots are associated to three-dimensional flows 
with topological entropy zero.

\begin{cor}
Given a nondegenerate gradient plane field flow on $S^3$, 
every component of the fixed point link is a zero-entropy knot.
\end{cor}
\pf 
Wada \cite{Wad89} classifies the knot types for 
cores of all round handle decompositions
on $S^3$. Each component is a zero-entropy knot.
\qed

\begin{rem}{\em
Much more can be said: Wada in fact classifies all possible {\em links}
which arise as round handle cores on all graph manifolds. This class
of {\em zero-entropy links} is an extremely restricted class, which lends
credence to the motto that simple dynamics implicate simple links in
dimension three. We note this same class of links appears independently 
in the study of nonsingular Morse-Smale flows, 
suspensions of zero-entropy disc maps, 
and in Bott-integrable Hamiltonian flows with two degrees of freedom. 
}\end{rem}

\begin{rem}{\em \label{rem_GP}
It is possible to construct gradient flows on $S^3$ (for example) 
in which the fixed point set is an embedded link which is {\em not}
a zero-entropy link. Let $L_1$ and $L_2$ denote any pair of links in $S^3$
which each have at least three components. Grayson and Pugh \cite{GP93} 
prove the existence of $C^\infty$ 
functions $\Psi_1, \Psi_2:\real^3\ra\real$ which 
have $L_1$ and $L_2$ as the (respective) sets of critical points. 
Moreover, these functions are proper and, for large enough $c\in\real$,
the inverse image of $c$ is a smooth 2-sphere near infinity. Hence, 
we may consider the balls $B_i$ bounded by $\Psi_i\inv(c)$ and glue
them together along the boundaries, obtaining $S^3$. The 
resulting function $\Psi$  given by $\Psi_1$ on $B_1$ and $-\Psi_2+2c$
on $B_2$ has as its gradient flow the split (unlinked) sum of
$L_1$ and $L_2$ as its fixed points. Thus, this flow cannot live 
within a plane field. 
}\end{rem}

It is not ostensibly clear that every zero-entropy link in $S^3$ is
realized as the zero set of a gradient flow within a plane field. 
We close this section with a realization theorem for such flows 
which shows that, in fact, a particular subclass of round-handle 
decompositions (and, hence, zero-entropy links) is realized.

\begin{lemma} \label{lem_Essential}
If $X$ is a nondegenerate gradient field on $M$ contained in 
the plane field $\eta$, then each index-1 round-handle $H$ in 
the decomposition must be attached to $\del M_i$ along annuli which are 
essential (homotopically nontrivial) in $\del M_i$.
\end{lemma} 
\pf
Assume that $M_{i}$ is the $i$th stage in a round handle 
decomposition, and that $H$ is an index-1 round handle 
with an exit annulus $E$ which is essential in $\del H$
by definition. By Theorem~\ref{thm_RH}, the intersection of 
$\eta$ with $\del M_i$ is always transverse. Thus, 
if $H$ is attached to $M_i$ along an annulus $A\subset\del M_i$, 
then the foliations given by the intersections of $\eta$ 
with the tangent planes to $A$ and $E$ respectively 
must match under the attachment. We claim this is 
impossible when $A$ is homotopically trivial in $\del M_i$.

Define the {\em index} of a smooth (oriented) curve $\gamma$ in 
an orientable surface with a (nonsingular, oriented) foliation $\F$  
to be the degree of the map which associates 
to each point $p\in\gamma$ the angle between  
the tangent vectors to $\gamma$ and $\F$ at $p$. This index is 
independent of the metric chosen and also 
invariant under homotopy of $\gamma$ or of $\F$; 
hence, we can speak of the index of an annulus in a surface
with foliation.

When $A$ is homotopically trivial, the index must be 
equal to $\pm 1$, since a foliation is locally a product. 
However, the index of the exit annulus $E\subset\del H$ must be 
zero as follows. 
Under the gradient field $X$, the core of the 1-handle is a curve
$\kappa$ of fixed points with transverse index 1 whose unstable manifold 
$W^u(\kappa)$ intersects $\del H$ transversally along the core of the 
exit set $E$. Deformation retract $E$ to a small neighborhood of 
$\W^u(\kappa)\cap\del H$ --- here, the intersections with $\eta$ 
are always transverse. Next, homotope the annulus to a neighborhood 
of $\kappa$ by integrating the gradient field $X$ backwards in 
time. This has the effect of taking the annulus transverse
to $W^u(\kappa)$ and sliding along $W^u(\kappa)$ back to $\kappa$.
Since $X$ points outwards along $W^u(\kappa)$, the image of the
annulus $E$ under the homotopy is always transverse to $X$, and
hence to $\eta$. The fact that $\eta\pitchfork\kappa$ then implies that 
the foliation on $E$ induced by $\eta$ must be homotopic through nonsingular
foliations to a product foliation by intervals on the annulus, 
which implies that the longitudinal annulus $E$ has index zero. Note that 
this works for exit sets $E$ which wrap any number of times about 
the longitude of $H$ (to cover both types of index-1 round handles).
\qed

Hence, any round-handle decomposition which is realizable 
as a gradient plane field flow must have all 1-handles attached 
along essential annuli. We call such a round-handle 
decomposition {\em essential}. 

\begin{thm}\label{thm_Realized}
Let $M$ be a compact 3-manifold with $L$ an indexed link. Then 
$L$ is realized as the indexed set of zeros for some nondegenerate 
gradient plane field flow on $M$ if and only if $L$ is the indexed 
set of cores for an essential RHD on $M$.
\end{thm}
\pf
The necessity is the content of Lemma~\ref{lem_Essential}. Given
any essential RHD, we construct a corresponding plane field gradient
flow. One may begin with the fact proved by Fomenko that
any essential round-handle decomposition can be generated by 
a vector field $X$ integrable via a Bott-Morse function 
$\Psi:M\ra\real$ with all critical sets being circles 
(see \cite{CMAN} for a 
detailed exposition). After choosing a metric on $M$ we claim that 
$-\grad\Psi$ lives within the plane field $\eta$ orthogonal to $X$.
Indeed, away from $L$ the plane field $\eta$ will be spanned by 
$\nabla\Psi$ and $\nabla\Psi\times X$ since these are linear 
independent vectors orthogonal to $X$ (recall $X$ is tangent
to the level sets of $\Psi$).  Thus $-\nabla\Psi$ clearly
lies in $\eta$ on the complement of $L$.  Along $L$ the
gradient $-\nabla\Psi$ lies in $\eta$ since it is zero.
\qed

The above construction may be modified so as to force the plane field to
twist monotonically along orbits of the gradient field, by a careful choice of the
vector field $X$. This implies that all of the permissible round handle decompositions 
may be realized by a totally nonintegrable plane field (a contact structure). Totally 
integrable plane fields are not so flexible, as will be illustrated in the next
section.

\section{Flows on foliations and contact structures}
\label{sec_5}

\subsection{Foliations}

In the case where our given plane field has some geometrical
property, we may further
restrict the types of round-handle decompositions which may 
contain a gradient flow. For example, if the plane field $\eta$
is integrable, it determines a foliation on the manifold. In this 
subsection, we note that, in this case, $S^3$ cannot support such 
a gradient flow. This result, which is an obvious corollary of
Novikov's Theorem on foliations, generalizes to other
three-manifolds.

\begin{figure}[htb]
\begin{center}
	\epsfxsize=3.45in\leavevmode\epsfbox{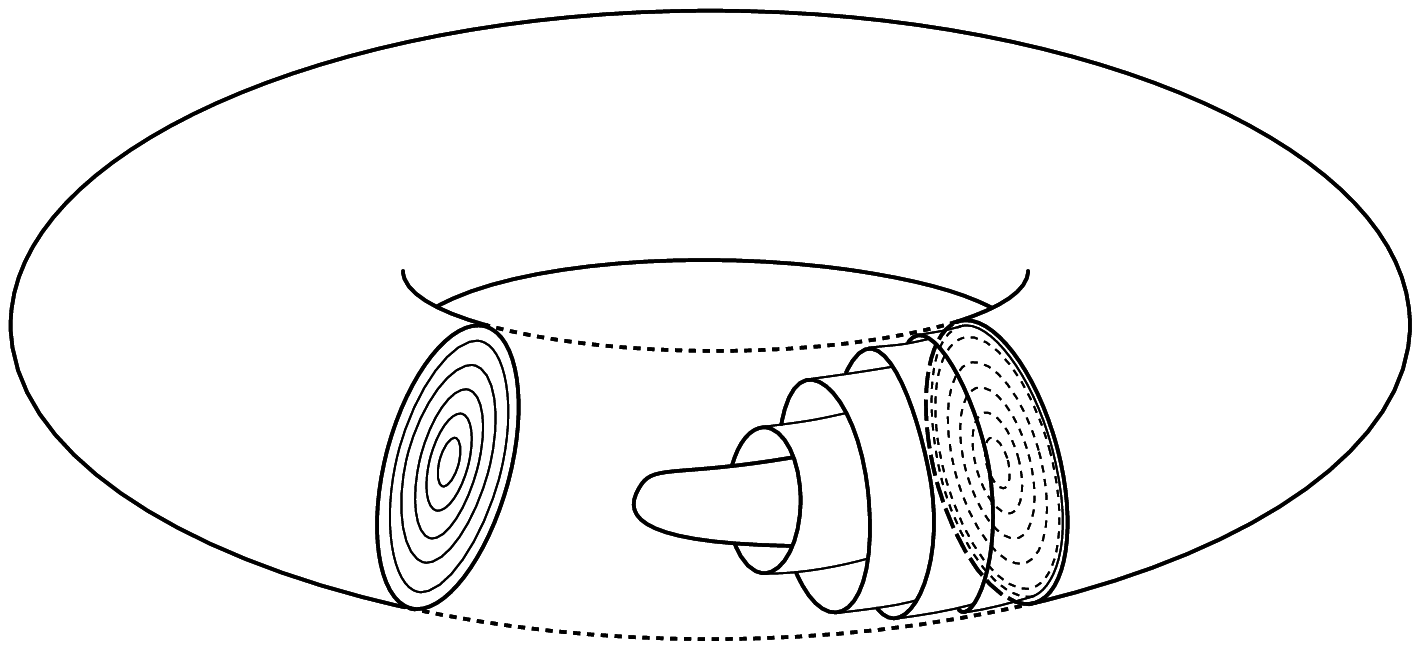}
\hspace{0.2in}
	\epsfxsize=1.15in\leavevmode\epsfbox{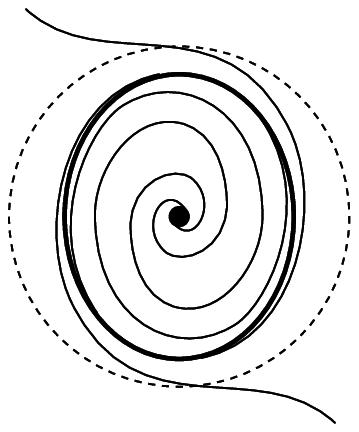}
\end{center}
\caption{A Reeb component in a foliation on a 3-manifold
	(left) can be perturbed into an overtwisted contact
	structure (right).} 
\label{fig_Reeb}
\end{figure}

Recall from the theory of foliations on three-manifolds 
(see, \eg, \cite{God91}) 
that a {\em Reeb component} is a foliation of the solid torus 
$D^2\times S^1$ that consists of the boundary $T^2$ leaf
along with a one-parameter family of leaves, each homeomorphic 
to $\real^2$ and limiting onto the boundary with nontrivial 
holonomy, as in Figure~\ref{fig_Reeb}. A codimension-one 
foliation of a three-manifold is {\em taut} if there do not
exist Reeb components or ``generalized'' Reeb components (see, 
\eg, \cite{ET97} for definitions). An equivalent definition 
of taut is that given any leaf $\L$ there exists a closed 
curve through $\L$ transverse to the foliation.

It is straightforward to show that gradient fields must lie within 
taut foliations:
\begin{thm}\label{thm_Taut}
Let $\cal F$ denote a codimension-1 foliation on a compact 
three-manifold $M$ which contains a nondegenerate gradient 
vector field $X$. Then $\F$ is taut.
\end{thm}
\pf
Assume that $X=-\grad\Psi$ is a nondegenerate gradient field on a foliation 
$\cal F$. For $\L$ a leaf of $\F$, the restriction of $X$ to $\L$ 
must also be a gradient flow. 
In the case where $\L$ is compact, there must be a nondegenerate 
fixed point of $X$ on $\L$ which lies on a circle of fixed points
transverse to $\F$ (note that in the case of a boundary torus in 
a Reeb component, this is an immediate contradiction). 
In the case where $\L$ is not a compact leaf, choose some nontrivial 
path $\gamma\subset\L$ whose endpoints are directly above one 
another in a local product chart. Then, by perturbing $\gamma$ 
to be transverse to $\F$, we may close it up to a transverse loop
through $\L$.
\qed

\begin{figure}[htb]
\begin{center}
	\epsfxsize=2.0in\leavevmode\epsfbox{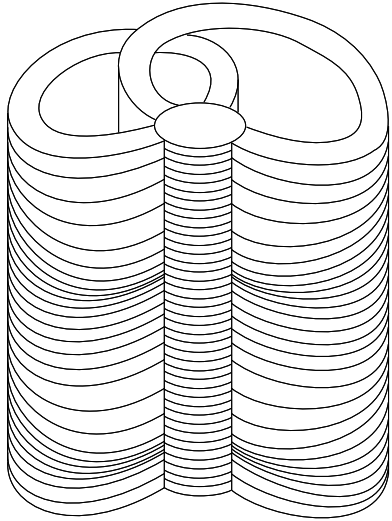}
\end{center}
\caption{A foliated RHD on $T^2\times S^1$ with holonomy along
each 1-handle: the 
2-handle has been removed and the $S^1$-factor cut open, revealing
a pair of 1-handles attached to a 0-handle with inverse attaching maps.}
\label{fig_Holonomy}
\end{figure}

This result can be greatly improved by considering the holonomy
of the foliation. Recall that the {\em holonomy} of any closed curve 
$\gamma:S^1\ra\L$ in a leaf $\L$ of a codimension-one 
foliation $\F$ is the germ of the Poincar\'e 
map associated to the characteristic foliation on an annulus 
transverse to $\L$ along $\gamma$.  
The holonomy of a curve is an invariant of its homotopy class within
the leaf. A foliation has {\em vanishing holonomy} if the holonomy of 
every curve $\gamma$ is trivial (the identity). 
\begin{thm}\label{thm_Foliation}
Any closed orientable three-manifold $M$ containing a nondegenerate
gradient field within a ($C^r$ for $r>1$) codimension-one foliation  
is a surface bundle over $S^1$.
\end{thm}
\pf
Suppose that $M$ admits a foliation $\F$ which supports
a nondegenerate gradient field. Then, by Theorem~\ref{thm_RH},
$M$ has a round handle decomposition where all the regular 
tori are transverse to the foliation. The foliation on each
round handle is equivalent to the product foliation by discs on 
$D^2\times S^1$, since these solid tori are filled with leaves 
transverse to the boundary each having a gradient flow with a single
fixed point. We show in subsequent steps that the foliation $\F$ 
may be modified within the round handle structure so that the new 
foliation $\F'$ has no holonomy. Once we show this, 
the celebrated theorem of Sacksteder implies that this foliation 
must be topologically conjugate to the kernel of a closed 
nondegenerate 1-form on $M$ \cite{Sac65}. The existence of this
1-form implies, via the theorem of Tischler 
\cite{Tis70}, that $M$ must be a surface bundle over $S^1$.
We illustrate in Figure~\ref{fig_Holonomy} below 
that it {\em is} possible to have gradient fields within a 
foliation having holonomy, so it is truly necessary to develop 
the following modification procedure, which makes use of 
``shearing'' the foliation along 1-handles (\cf \cite{EHN81}).

Denote by $M_0$ the (disjoint) union of all the 0-handles and 
by $M_i$ ($1\leq i\leq N$) the subsequent stages in the decomposition:
\[
	M_i = \left(\left. M_{i-1}\sqcup H_i\right) \right/ \phi_i,
\]
where $H_i$ is the $i$th 1-handle and $\phi_i:E_i\hra\del M_{i-1}$ 
is the attaching map on the exit set $E_i\subset H_i$.
Recall that each exit set $E_i$ is either one or two annuli and  
that the boundary of each $M_i$, $\del M_i$, is the disjoint union 
of a collection of tori. 

For each 1-handle $H_i$, let $W_i^u$ denote the (2-dimensional) 
unstable manifold to the core of $H_i$. Modify the 
round handle structure so that each $H_i$ is very
``thin'' -- that is, each $H_i$ is restricted to a small neighborhood of 
$W^u_i$, appending the ``leftover'' portion to the neighboring 2-handles.
Denote by 
\[
	{\mbox{Bd}}_i = \del M_0\bigcup_{j=1}^i W^u_j ,	
\]
the 2-complex given by the union of all the 0-handle boundaries and 
unstable manifolds of the 1-handles in $M_i$. 

{\em Claim 1:} 
$\F$ has vanishing holonomy if the restriction of $\F$ to Bd$_N$ 
has vanishing holonomy.

{\em Proof 1:} Let $\gamma$ denote a loop within a leaf 
$\L$ of $\F$. Then the restriction of $\L$ to each $k$-handle 
is a collection of disjoint discs whose boundaries lie in the union of 
0- and 1-handles. Push $\gamma$ to these boundaries and, since 
$\del H_i$ is very close to $W^u_i$, perturb $\gamma$ to 
lie within $W^u_i$ for each $H_i$ it intersects.  
Since Bd$_N$ is transverse to $\F$, we may choose the transverse 
annulus $A$ containing $\gamma$ to lie within this set.
\qed$_1$

In what follows, we consider holonomy on the 2-complex 
Bd$_N$, keeping in mind that the 1-handles are actually thin 
neighborhoods of the 2-cells $W^u_i$. The holonomy on each 
component of $\del M_i$ is equivalent to that on the corresponding
piece of Bd$_N$ since each $H_i$ has a product foliation.

{\em Claim 2:} The maps $\{\phi_i\}_1^N$ may be isotoped so that 
the induced foliation $\F'$ on Bd$_N$ is without holonomy.

{\em Proof 2:}
It suffices to show that the foliation restricted to each 
$\del M_i$ is without holonomy (a product foliation): we
proceed by induction on $i$. 
On the boundary of $M_0$ the foliation $\F$ restricts to a product 
foliation by circles. 
Assume as an induction hypothesis a lack of holonomy on $\del M_{i-1}$. 
There are three cases to consider: (1) $H_i$ is an orientable handle 
with attaching circles in the same component of $\del M_{i-1}$;
(2) $H_i$ is orientable with attaching circles in two 
distinct components of $\del M_{i-1}$; and (3) $H_i$ is 
nonorientable. 

Case (1): Let $C_\pm$ denote the circles in the selected
component $T$ of $\del M_{i-1}$ along which $W^u_i$ is attached. 
Note $C_\pm$ divides $T$ into two annuli $A_0$ and $A_1.$ 
After fixing a diffeomorphism from $C_+$ to $C_-$ there is a 
``handle holonomy map'' 
$f_H:C_+\ra C_-$ which is the diffeomorphism given by sliding 
along leaves on $\phi_i(H_i)$. There are corresponding 
``boundary holonomy maps'' $f_j:C_+\ra C_-$ given by sliding
along leaves on $A_j$. Isotope $\phi_i$ on $C_+$ so that 
$f_H$ equals $f_0$ up to a rigid rotation 
(which is necessary in order to add subsequent handles along curves
transverse to $\F$ --- see Claim 3). 
The holonomy on the two new components of  $\del M_i$ is determined 
by taking the transverse 
curve $C_+$ (actually one must take a parallel copy of $C_+$ that 
sits in $\del M_i$) as a section. These holonomy maps factor as 
$f_1\inv\circ f_H$ and $f_H\inv\circ f_0$; 
however, the holonomy along $C_+$ within $\del M_{i-1}$ is a map of the 
form $f_1\inv\circ f_0$, which, by induction, is a rigid rotation. 
Hence, up to rotations, $f_0=f_1$. Since we chose $f_H=f_0$ 
up to rotations the holonomy on $\del M_i$ vanishes.

Case (2): If $H_i$ connects two disconnected boundary components
of $M_{i-1}$, then the holonomy along $H_i$ will always cancel 
with itself as follows. Denote by $f_H:C_+\ra C_-$ the handle holonomy 
maps as before. Then the global holonomy map along $\del M_i$ 
is of the form $g_+\circ f_H\circ g_- \circ f_H\inv$, 
where $g_+:C_+\ra C_+$ and $g_-:C_-\ra C_-$ are holonomy self-maps 
along loops in the two components of $\del M_{i-1}$, and hence by 
induction, identity maps. 

Case (3): If $H_i$ has connected exit set, the proof follows as 
in Case (1), since the handle must connect a single component
of $\del M_{i-1}$ to itself: isotope $\phi_i$ so that the handle holonomy 
map equals the holonomy map along the boundary up to a rigid 
rotation. 
\qed$_2$

{\em Claim 3:}
This ``linearization'' of $\F$ does not affect the topology of $M$. 

{\em Proof 3:}
Throughout the addition of the
1-handles, nothing about the topology of $M$ has changed, since
the handle structure is identical --- we modify only the foliation.
However, after attaching the last 1-handle, the characteristic 
foliation on the boundary tori must be linear and rational, in order
to glue in the 2-handles respecting the product foliation on their 
boundaries. The slopes of $\F$ restricted to $\del M_N$ completely
determine the topology of $M$ after adding the 2-handles (these 
are Dehn filling coefficients).

Hence, we must be able to linearize all of the attaching maps for the 
1-handles without changing the boundary slopes at the end of the
sequence. To do so, we preserve at every stage the {\em rotation 
number} of the holonomy maps $h_i$ which slide the attaching 
curves of $H_i$ along $\del M_i$. Recall that to every 
diffeomorphism $f:S^1\ra S^1$ is 
associated a rotation number $\rho_f\in\real/\zed$ which measures
the average displacement of orbits of $f$ (see, \eg, \cite{GH83}). 
When modifying $\phi_i$ to $\tilde\phi_i$ in the above
procedure, we may compose $\tilde\phi_i$ with a rigid rotation 
by the angle necessary to preserve the rotation number of the 
holonomy map $h_i$ acting on the attaching curves in $\del M_{i-1}$
(without adding further Dehn twists). 
This shearing maintains the average slope of the boundary foliation 
at each stage without adding holonomy. Hence, at the end 
of the 1-handle additions, when the original foliation had all 
boundary components with linear foliations of a particular fixed 
slope, the modified foliation also has linear boundary foliations 
with the same slope. Thus, adding the 2-handles is done using the 
same surgery coefficients, yielding the original manifold $M$ 
with a foliation having trivial holonomy. 
\qed$_3$

Claims 1-3 complete the proof of Theorem~\ref{thm_Foliation}.
\qed

\begin{rem}{\em 
Of course, not every surface bundle over $S^1$ may support a 
gradient field within a foliation: there is still the restriction that
$M$ be a graph-manifold. This translates precisely into a condition 
on the monodromy map of the fibration --- the monodromy must
be of periodic (or reducibly periodic) type with respect to the 
Nielsen-Thurston classification of surface homeomorphisms. 
Any pseudo-Anosov piece in the monodromy forces
hyperbolicity, contradicting the graph condition. 
It is not hard to see that any such bundle can be 
given a gradient field lying within each fiber $F$ of the bundle by 
choosing a Morse function $\phi:F\ra\real$ which is equivariant with 
respect to the monodromy map.
}\end{rem}

\begin{rem}{\em
All of the results of this section apply not only to 
gradient flows, but also to {\em gradient-like flows}, or
flows for which there exists a function which decreases strictly 
along non-constant flowlines. The reason why nondegenerate gradient-like
flows in foliations determine round-handle decompositions whereas 
for general plane fields they do not lies in the fact that 
the Hopf bifurcation of Proposition~\ref{prop_Hopf} cannot take 
place among gradient-like flows in the integrable case, while 
it can in the nonintegrable.
}\end{rem}

\subsection{Contact structures}

In contrast to the case of an integrable plane field, one may 
consider the class of {\em contact structures}, which has 
attracted interest in the fields of symplectic geometry and
topology, knot theory, mechanics, and hydrodynamics.
\begin{dfn}\label{def_Contact}{\em 
A {\em contact form} on a three-manifold $M$ is a one-form 
$\alpha\in\Omega^1(M)$ such that the Frobenius integrability 
condition fails everywhere: that is,
\be 
	\alpha\wedge d\alpha \neq 0.
\ee
A {\em contact structure} on $M$ is a plane field $\xi$ which is 
the kernel of a locally defined contact form: that is,
\be
	\xi_p = \{ \bfv \in T_p : \alpha(\bfv)=0 \},
\ee
for each $p\in M$. 
}\end{dfn}
Contact structures are thus maximally nonintegrable: the plane 
field is locally twisted everywhere. 
One may think of a contact structure as being an {\em anti-foliation},
which leads one to suspect that the topology of the manifold 
may be connected to the geometry of the structure, as is often
the case with foliations. Indeed, the contrast between 
foliations with Reeb components and those without Reeb components
is reflected in the {\em tight / overtwisted dichotomy} in contact
geometry (due primarily to Eliashberg \cite{Eli92} and Bennequin \cite{Ben82}). 

\begin{dfn}\label{def_Tight}{\em
Given a contact structure $\xi$ on $M$ and an embedded 
surface $F\subset M$, then the {\em characteristic foliation} 
$F_\xi$ is the (singular) foliation induced by the (singular) line field
$\{T_pF\cap\xi_p : p\in F\}$. 
A contact structure $\xi$ is {\em overtwisted} if there exists 
an embedded disc $D\in M$ such that $D_\xi$ has a limit cycle, 
as in Figure~\ref{fig_Reeb}(right).
A contact structure is {\em tight} if it is not overtwisted.
}\end{dfn}

The classification of contact structures follows along lines similar
to that of codimension one foliations with or without Reeb
components.  An infinite number of homotopically distinct overtwisted 
contact structures exist on every closed orientable 
three-manifold \cite{Mar71,Lut77} and are algebraically classified up to
homotopy \cite{Eli89}. Tight structures, on the other hand, are 
quite mysterious: \eg, it is unknown whether they exist on all 
three-manifolds. 

Several examples of the similarity between tight contact structures and
Reebless foliations are provided by the recent work of Eliashberg
and Thurston \cite{ET97}. For example, Reebless foliations can be
perturbed into tight contact structures and foliations with Reeb 
components can be perturbed into overtwisted structures (\cf 
Figure~\ref{fig_Reeb}). Also, both Reebless foliations 
and tight structures satisfy a strong inequality restricting Euler 
classes. Tight structures are somewhat more general than 
Reebless foliations since the former can exist on $S^3$ \cite{Ben82}
while the latter cannot \cite{Nov67}. Likewise, overtwisted 
structures are slightly more general than their foliation 
counterparts via the following observation, to be contrasted with 
Theorem~\ref{thm_Taut}:

\begin{prop}
Any nondegenerate gradient field $X$ which lies within a tight 
contact structure $\xi$ on $M^3$ also lies within an overtwisted
contact structure $\xi'$ on $M^3$.
\end{prop}
\pf
The canonical way to turn a tight structure into an 
overtwisted structure is by performing a {\em Lutz twist}
\cite{Lut77,Mar71} on a simple closed curve $\gamma$ 
transverse to $\xi$. We execute a version of this twisting which 
respects a gradient field.

Given a gradient field $X\subset\xi$, choose a curve $\gamma$ of fixed
points of index zero (sinks). Translate 
the function $\Psi$ whose gradient defines $X$ so that 
$\rest{\Psi}{\gamma}\equiv 0$. Since $\gamma$ is an index 
zero curve, $\Psi$ increases as one moves radially away from $\gamma$.

Let $N$ denote a tubular neighborhood of $\gamma$ whose boundary 
is a connected component of $\Psi\inv(\eps)$ for some $\eps>0$.
Place upon $N$ the natural cylindrical coordinates $(\Psi,\theta,z)$.
In analogy with Lemma~\ref{lem_Coords}, we may choose $\theta$ and 
$z$ so that $\rest{\xi}{N}$ is the kernel of the locally defined 1-form 
\be
	\alpha = g(\Psi,\theta,z)d\theta + dz ,
\ee
for some function $g$ with $g(0,\theta,z)=0$. 
The contact condition implies that 
\be\label{eq_grow}
	\frac{\del g}{\del\Psi} > 0 .
\ee
Replacing this structure locally with the kernel of the form 
\be 
\alpha' = 
	\sin\left(\frac{\pi}{4}+\frac{2\pi g}{g(\eps,\theta,z)}\right)
		g\,d\theta 
+ 	\cos\left(\frac{\pi}{4}+\frac{2\pi g}{g(\eps,\theta,z)}\right)dz ,
\ee
yields a contact structure since 
\be
\alpha'\wedge d\alpha' = \left[
	\cos\left(\frac{\pi}{4}+\frac{2\pi g}{g(\eps,\theta,z)}\right)
	\sin\left(\frac{\pi}{4}+\frac{2\pi g}{g(\eps,\theta,z)}\right)
	+ \frac{2\pi g}{g(\eps,\theta,z)}
\right] \frac{\del g}{\del\Psi}\,d\Psi\wedge d\theta\wedge dz ,
\ee
and this coefficient is positive by Equation~\ref{eq_grow}.
This contact structure agrees with that defined by $\alpha$ 
along the torus $\Psi=\eps$ since 
\be
	\rest{\alpha'}{\Psi=\eps} 
= 	\sin\left(\frac{9\pi}{4}\right)g(\eps,\theta,z)d\theta 
+ 	\cos\left(\frac{9\pi}{4}\right) dz 
=	\frac{\sqrt{2}}{2}\rest{\alpha}{\Psi=\eps} ,
\ee
and these have the same kernel.
Furthermore, this modified structure contains the
vector field $X=\grad\Psi$, since $X$ points in the $-d/d\Psi$ 
direction. Finally, one can easily show that a perturbation of a 
constant-$z$ disc has a limit cycle in the characteristic foliation 
near $c=\eps/2$ (\cf \cite{Ben82}); hence, this 
defines an overtwisted structure containing $X$. This construction 
can obviously be done in the $C^\infty$ category using bump
functions. 
\qed

\begin{figure}[htb]
\begin{center}
	\epsfxsize=3.0in\leavevmode\epsfbox{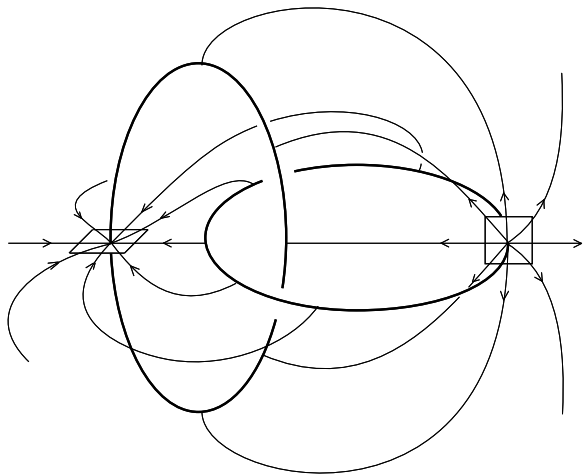}
\end{center}
\caption{The gradient field on $S^3$ having a Hopf link of fixed points
	exists within both tight and overtwisted contact structures.} 
\label{fig_NS}
\end{figure}

\begin{ex}{\em
Consider the flow on $S^3$ (considered as the unit sphere in
$\real^4$ with the induced metric) given by the gradient of the
function
\be
	\Psi = \frac{1}{2}(x_1^2+x_2^2)-\frac{1}{2}(x_3^2+x_4^2),
\ee
the gradient being taken in $S^3$. 
One can check that the fixed point set consists of a
pair of unknots linked once in a {\em Hopf link}, as in 
Figure~\ref{fig_NS}. The standard tight contact form on $S^3$ is
\be
	\alpha = \frac{1}{2}\left( 
	x_1 dx_2 - x_2 dx_1 + x_3 dx_4 - x_4 dx_3 \right) .
\ee
A simple calculation shows that $\alpha$ is a contact form on $S^3$
with $\grad\Psi\subset\ker\alpha$. However, we may Lutz twist this
structure in a neighborhood of the fixed point links: a family of
such overtwisted forms ($n\in\zed^+$) is given by \cite{GV83}
\be 
	\alpha_n = \cos\left(\frac{\pi}{4}+n\pi(x_3^2+x_4^2)\right)
	(x_1 dx_2 - x_2 dx_1) + 
	\sin\left(\frac{\pi}{4}+n\pi(x_3^2+x_4^2)\right)
	(x_3 dx_4 - x_4 dx_3), 
\ee
from which it can be shown that $\grad\Psi\subset\ker\alpha_n$.
Here, the integer $n$ denotes the number of twists that the
plane field undergoes as an orbit travels from source to sink
in Figure~\ref{fig_NS}.
}\end{ex}

\section{Two questions}
\label{sec_6}

This work has focused on the case of gradient flows in plane fields
in dimension three, as the round-handle theory is most interesting here. 
However, there are natural questions about gradient flows in 
arbitrary distributions for manifolds of dimension greater than 
three. We do not present any results in this area, but rather note
that many of the tools remain valid: fixed point sets of a vector
field constrained to a codimension-$k$ distribution consists of a 
finite collection of embedded $k$-dimensional submanifolds.

Two problems emerge. In the case of a codimension-one distribution, 
nondegenerate gradient fields induce round-handle decompositions.
However, every manifold of dimension greater than three whose Euler 
characteristic is zero possesses
an RHD. Are there any such manifolds of dimension greater than three
which do not possess a nondegenerate gradient field tangent to a 
codimension-one distribution? Secondly, in the case of 
higher codimension distributions, what restrictions exist on 
the topology of the fixed point sets? The case of a plane field on 
a four-manifold is particularly interesting with respect to the 
genera of the (two-dimensional) fixed point sets.


\noindent{\sc ACKNOWLEDGMENTS}
\vspace{0.05in}

This work has been supported in part by the National Science
Foundation [JE: grant DMS-9705949; RG: grant DMS-9508846]. 
The authors wish to thank Mark Brittenham, John Franks, Will Kazez, 
Alec Norton, and Todd Young for their input. Special thanks are due
the referee for constructive remarks.



\end{document}